\documentclass[12pt]{amsart}
\usepackage{mathrsfs}
\usepackage{txfonts}
\usepackage{amssymb}
\usepackage{color}
\makeatletter \@mparswitchfalse \makeatother
\newtheorem{theorem}{Theorem}[section]

\def\eqref#1{(\ref{eq#1})}

\numberwithin{equation}{section}
\begin{document}

\title[]{A counterexample to Arveson's finiteness question
for subdiagonal algebras}
\date{}
\author[ ]{Guoxing Ji }
  \address{School  of Mathematics and Statistics,
  Shaanxi Normal University,
  Xian , 710119, People's  Republic of  China}
  \email{gxji@snnu.edu.cn}

\begin{abstract}
   We construct a unital \(\sigma\)-weakly closed subalgebra \(\mathfrak A\) of a von Neumann factor \(\mathcal M\) of type \(\mathrm{II}_1\), whose diagonal is \(\mathfrak D=\mathfrak A\cap\mathfrak A^*=\mathbb CI\). We show that there exist uncountably many faithful normal conditional expectations from \(\mathcal M\) onto \(\mathfrak D\) with respect to each of which \(\mathfrak A\) is a maximal subdiagonal algebra. Let \(\tau\) denote the unique normalized faithful normal trace on \(\mathcal M\). The conditional expectation associated with \(\tau\) preserves this trace, whereas none of the other conditional expectations is preserved by any nonzero finite trace on \(\mathcal M\). This gives a negative answer to  Arveson's finiteness  question for subdiagonal algebras    in \emph{Amer. J. Math.} \textbf{89}(1967), 578--642.
\end{abstract}

\maketitle

\section{Introduction}
In \cite{arv1}, Arveson introduced the notion of subdiagonal algebras in the study of analyticity in operator algebras. Let $\mathcal M$ be a $\sigma$-finite von Neumann algebra acting on a complex Hilbert space $\mathcal H$, and let $\mathfrak A\subseteq\mathcal M$ be a $\sigma$-weakly closed unital subalgebra of $\mathcal M$. Set
\[
\mathfrak D=\mathfrak A\cap\mathfrak A^*.
\]
Suppose that there exists a faithful normal conditional expectation $\Phi:\mathcal M\to\mathfrak D$ such that the following conditions hold:

(i) $\Phi$ is a homomorphism on $\mathfrak A$;

(ii) $\mathfrak A+\mathfrak A^*$ is $\sigma$-weakly dense in $\mathcal M$.

Then $\mathfrak A$ is called a subdiagonal algebra of $\mathcal M$ with respect to $\Phi$.

We say that $\mathfrak A$ is a maximal subdiagonal algebra of $\mathcal M$ with respect to $\Phi$ if there is no strictly larger subdiagonal algebra of $\mathcal M$ with respect to the same conditional expectation $\Phi$ that contains $\mathfrak A$.

Put
\[
\mathfrak A_0=\{X\in\mathfrak A:\Phi(X)=0\}
\]
and
\[
\mathfrak A_m=
\left\{
X\in\mathcal M:
\Phi(AXB)=\Phi(BXA)=0,\
\forall A\in\mathfrak A,\ \forall B\in\mathfrak A_0
\right\}.
\]
By \cite[Theorem 2.2.1]{arv1}, $\mathfrak A_m$ is the unique maximal subdiagonal algebra of $\mathcal M$ with respect to $\Phi$ containing $\mathfrak A$.

Throughout this paper, the term subdiagonal algebra will always refer to a subdiagonal algebra with respect to a faithful normal conditional expectation.

Following Arveson \cite[Definition 4.1.0]{arv1},
a subdiagonal algebra $\mathfrak A$ with respect to $\Phi$ is called finite if there exists
a faithful normal finite trace $\tau$ on $\mathcal M$ such that
$\tau\circ\Phi=\tau$.  Exel proved that a finite subdiagonal algebra  is automatically maximal in
\cite[Theorem 7]{Exel}.  Arveson asked whether every subdiagonal
algebra in a finite von Neumann algebra is finite; see
\cite[Remark 4.1.1]{arv1}.

There are positive answers under some additional
assumptions. Kawamura and Tomiyama proved that if $\mathcal M$
is a finite von Neumann algebra acting on a separable Hilbert
space, $\alpha$ is a $\sigma$-weakly continuous flow on
$\mathcal M$, and the center $Z(\mathcal M)$ is $\alpha$-finite,
then the analytic operator algebra $H^\infty(\alpha)$ is a finite maximal
subdiagonal algebra; see \cite[Proposition 2.7]{KawamuraTomiyama}.
Here $\alpha$-finiteness means that the family of
$\alpha$-invariant normal states separates the positive elements
of the center.  On the other hand, It is proved that
  every type $1$ subdiagonal
algebra in a finite von Neumann algebra is finite   in
\cite[Theorem 3.1]{Ji}.

In this note, we give a counterexample to Arveson's finiteness problem.
The example is constructed directly on $L^2(\mathbb T^2)$.
Its von Neumann algebra is unitarily equivalent to the GNS
weak closure of the irrational rotation $C^*$-algebra, which
is the hyperfinite factor of type $\mathrm{II}_1$.

\section{A counterexample}

Fix an irrational number $\theta$ and put
\(
\lambda=e^{2\pi i\theta}.
\)
Let $m$ be the normalized Haar measure on $\mathbb T$, and
consider $L^2(\mathbb T^2,m\times m)$. Define two operators
$U,V\in \mathcal B(L^2(\mathbb T^2))$ by
\[
(Uf)(z,w)=zf(z,w),\qquad
(Vf)(z,w)=wf(\lambda^{-1}z,w),
\quad f\in L^2(\mathbb T^2).
\]
Both $U$ and $V$ are unitary, and a direct computation gives
\begin{equation}\label{f1}
UV=\lambda VU.
\end{equation}
Let
\[
\mathcal M=W^*(U,V)\subseteq \mathcal B(L^2(\mathbb T^2)),
\]
and denote by $1$ the constant function on $\mathbb T^2$.

For all $n,m\in\mathbb Z$,
\[
U^nV^m1=z^nw^m.
\]
Hence the vectors $U^nV^m1$, $n,m\in\mathbb Z$,
form the standard orthonormal basis of $L^2(\mathbb T^2)$.
In particular, $1$ is cyclic for $\mathcal M$.

Define
\begin{equation}\label{f2}
\tau(X)=
\int_{\mathbb T^2}(X1)(z,w)\,dm(z)\,d m(w)=\langle X1,1\rangle,
\qquad X\in\mathcal M,
\end{equation}
where $\langle \cdot, \cdot\rangle $ denotes the inner product of $L^2(\mathbb T^2)$. Then
\begin{equation}\label{f3}
\tau(U^nV^m)=
\begin{cases}
1,&n=m=0,\\
0,&\text{otherwise}.
\end{cases}
\end{equation}
Relation (\ref{f1}) yields
\begin{equation}\label{f4}
(U^nV^m)(U^kV^\ell)
=\lambda^{-mk}U^{n+k}V^{m+\ell},
\qquad n,m,k,\ell\in\mathbb Z.
\end{equation}
It follows from   (\ref{f3}) and (\ref{f4}) that $\tau$ is tracial on the
algebraic $*$-algebra generated by $U$ and $V$. Since $\tau$
is a vector state, it is normal. The algebraic $*$-algebra
generated by $U,V$ is $\sigma$-weakly dense in $\mathcal M$,
so separate $\sigma$-weak continuity of multiplication implies
that $\tau$ is tracial on all of $\mathcal M$.

We identify the factor in this construction by a standard
reference. Let $A_\theta$ be the irrational rotation
$C^*$-algebra generated by unitaries $\widetilde U,\widetilde V$
satisfying
\[
\widetilde U\widetilde V
=\lambda\widetilde V\widetilde U,
\]
and let $\tau_\theta$ be its canonical faithful normalized
trace. The unitary which sends
\[
z^nw^m\longmapsto
\pi_{\tau_\theta}(\widetilde U^n\widetilde V^m)
\Omega_{\tau_\theta}
\]
intertwines the operators $U,V$ above with the left GNS
representations of $\widetilde U,\widetilde V$, respectively.
Thus $\mathcal M$ is unitarily equivalent to
$\pi_{\tau_\theta}(A_\theta)''$. Goswami and Sahu identify
this GNS weak closure as the hyperfinite factor of type
$\mathrm{II}_1$; see \cite[Section 2.2]{GoswamiSahu}.
The construction of $A_\theta$ and its canonical trace is
also discussed in \cite{Rieffel}.
For the general theory of quantum tori, we refer the reader to \cite{ChenXuYin}. Consequently, $\tau$ is
the faithful normal tracial state on $\mathcal M$. Since $\mathcal M$ is a factor, the normalized faithful normal trace on $\mathcal M$ is unique.

Put
\[
P=\{(n,m)\in\mathbb Z^2:n>0\}
\cup\{(0,m)\in\mathbb Z^2:m\geq0\},
\]
and define
\begin{equation}\label{f5}
\mathfrak A=
\overline{\operatorname{span}}^{\,\sigma\text{-weak}}
\{U^nV^m:(n,m)\in P\}.
\end{equation}
Since $P+P\subseteq P$, equation (\ref{f4}) shows that
$\mathfrak A$ is a unital $\sigma$-weakly closed subalgebra
of $\mathcal M$.

For $X\in\mathcal M$, define its Fourier coefficients by
\begin{equation}\label{f6}
\widehat X(n,m)=\tau\big((U^nV^m)^*X\big),
\qquad n,m\in\mathbb Z.
\end{equation}
The monomials $U^nV^m$ form an orthonormal basis of
$L^2(\mathcal M,\tau)$, and each coefficient functional in
(\ref{f6}) is normal. It follows from (\ref{f5}) that
\[
X\in\mathfrak A
\quad\Longrightarrow\quad
\widehat X(n,m)=0\quad\text{if }(n,m)\notin P.
\]
Similarly, if $X\in\mathfrak A^*$, then
\[
\widehat X(n,m)=0\quad\text{if }(n,m)\notin -P.
\]
Since
\[
P\cap(-P)=\{(0,0)\},\qquad P\cup(-P)=\mathbb Z^2,
\]
we obtain
\begin{equation}
\label{f7}
\mathfrak D=\mathfrak A\cap\mathfrak A^*=\mathbb CI.
\end{equation}
Moreover, every monomial $U^nV^m$ belongs to either
$\mathfrak A$ or $\mathfrak A^*$. Their linear span is
$\sigma$-weakly dense in $\mathcal M$, and hence
\begin{equation}\label{f8}
\overline{\mathfrak A+\mathfrak A^*}^{\,\sigma\text{-weak}}
=\mathcal M.
\end{equation}
Formula (\ref{f3}) implies that  $\tau$ is multiplicative on $\mathfrak A$.  Thus if we define $\Phi_0(A)=\tau(A)I$ for all $A\in\mathcal M$, then  $\Phi_0$ is a faithful normal conditional expectation from  $\mathcal M$ onto the diagonal $\mathfrak D=\mathbb CI$ of $\mathfrak A$,  which is a homomorphism on $\mathfrak A$. It follows that $\mathfrak A$ is a finite subdiagonal algebra with respect to $\Phi_0$.

Next we construct  a faithful normal conditional  expectation $\Phi_r$  for each $0< r<1$,  such that $\mathfrak A$ is also a subdiagonal algebra with respect to $\Phi_r$ which is not finite.  We define
\begin{equation}\label{f9}
H_r=(1-r^2)(I-rV)^{-1}(I-rV^*)^{-1}.
\end{equation}
The functional calculus for $V$ gives the norm-convergent
series
\begin{equation}\label{f10}
H_r=\sum_{j\in\mathbb Z}r^{|j|}V^j.
\end{equation}
  Moreover,
\begin{equation}\label{f11}
\frac{1-r}{1+r}I\leq H_r\leq
\frac{1+r}{1-r}I,\qquad \tau(H_r)=1.
\end{equation}
Set
\begin{equation}\label{f12}
\varphi_r(X)=\tau(H_rX),\qquad X\in\mathcal M.
\end{equation}
For every $X\in\mathcal M_+$, traciality and (\ref{f11}) imply
\[
\varphi_r(X)
=\tau(X^{1/2}H_rX^{1/2})
\geq\frac{1-r}{1+r}\tau(X).
\]
Thus $\varphi_r$ is a faithful normal state on $\mathcal M$.
Since $\mathfrak D=\mathbb CI$, the map
\begin{equation}\label{f13}
\Phi_r(X)=\varphi_r(X)I,\qquad X\in\mathcal M,
\end{equation}
is a faithful normal conditional expectation from $\mathcal M$
onto $\mathfrak D$.

By (\ref{f3}) and (\ref{f10}),
\begin{equation}\label{f14}
\varphi_r(U^nV^m)=
\begin{cases}
r^{|m|},&n=0,\\
0,&n\neq0.
\end{cases}
\end{equation}
We next verify that $\varphi_r$ is multiplicative on
$\mathfrak A$. Let $(n,m),(k,\ell)\in P$. If $n+k>0$,
then at least one of $n,k$ is positive. Equations  (\ref{f4}) and  (\ref{f14}) give
\[
\varphi_r\big((U^nV^m)(U^kV^\ell)\big)
=\varphi_r(U^nV^m)\varphi_r(U^kV^\ell)=0.
\]
If $n+k=0$, then $n=k=0$ and $m,\ell\geq0$, whence
\[
\varphi_r(V^{m+\ell})
=r^{m+\ell}
=r^m r^\ell
=\varphi_r(V^m)\varphi_r(V^\ell).
\]
Therefore multiplicativity holds on the linear span in (\ref{f5}).
Normality of $\varphi_r$ and separate $\sigma$-weak
continuity of multiplication extend it to all of
$\mathfrak A$:
\begin{equation}\label{f15}
\varphi_r(XY)=\varphi_r(X)\varphi_r(Y),
\qquad X,Y\in\mathfrak A.
\end{equation}
Consequently,
\begin{equation}\label{f16}
\Phi_r(XY)=\Phi_r(X)\Phi_r(Y),
\qquad X,Y\in\mathfrak A.
\end{equation}
Equations (\ref{f7}), (\ref{f8}), and (\ref{f16}) show that $\mathfrak A$ is a
subdiagonal algebra with respect to every $\Phi_r$.  It is known that $\mathfrak A$ is maximal with respect to $\Phi_r$(cf. \cite[p. 138]{Xu}). Since $\tau$ is the unique normalized faithful normal trace on \(\mathcal M\) and $\tau\circ\Phi_r\ne\tau$, $\mathfrak A$ is not finite.

\begin{theorem}\label{t1}
For every $0< r<1$, the algebra $\mathfrak A$ is  a maximal
subdiagonal algebra with respect to $\Phi_r$ which is  not finite.
\end{theorem}

  Theorem \ref{t1}  gives a negative answer  to Arveson's finiteness problem in  \cite[Remark 4.1.1]{arv1}.

\noindent{\bf Acknowledgment}  The development and refinement of the mathematical arguments in this manuscript benefited from discussions with ChatGPT (GPT-6). All AI-assisted contributions have been carefully reviewed and verified by the authors, who take full responsibility for the content of the manuscript.

\end{document}